\documentclass[10pt]{article}
\usepackage{amssymb, amsmath, amsthm, color,graphicx, amsfonts,a4}
\usepackage{amscd}

\definecolor{Red}{cmyk}{0,1,1,0}

\definecolor{verde}{cmyk}{1,0,1,0}

\definecolor{azul}{cmyk}{1,1,0,0}

\numberwithin{equation}{section}

\def\cal{\mathcal}

\newcommand{\lra}{\leftrightarrow}

\DeclareMathOperator{\sen}{sen}

\newcommand{\E}{\mathbb{E}}
\newcommand{\N}{\mathbb{N}}

\newcommand{\Z}{\mathbb{Z}}
\renewcommand{\P}{\mathbb{P}}

\newcommand{\V}{\mathbb{V}}

\newcommand{\be}{\begin{equation}}
\newcommand{\ee}{\end{equation}}

\newtheorem{teorema}{Theorem}
\newtheorem{proposicao}{Proposition}

\newtheorem{lema}{Lemma}

\begin{document}
\title{Anisotropic Percolation on Slabs}
\author{Rodrigo G. Couto \footnote{ Departamento de
Matem{\'a}tica, Universidade Federal de Minas Gerais, Av.
Ant\^onio Carlos 6627 C.P. 702 CEP30123-970 Belo Horizonte-MG,
Brazil and Departamento de Matem\'atica, Universidade Federal de Ouro Preto, Rua Diogo de Vasconcelos 122 CEP35400-000 Ouro Preto-MG} ,
Bernardo N. B. de Lima \footnote{ Departamento de Matem{\'a}tica, Universidade Federal de Minas Gerais, Av. Ant\^onio
Carlos 6627 C.P. 702 CEP30123-970 Belo Horizonte-MG, Brazil} ,
R\'emy Sanchis$\;^{\dag}$}
\maketitle
\begin{abstract}
We consider anisotropic
independent bond percolation models on the slab $\Z^2\times\{0,\dots,k\}$, where we suppose that the axial (vertical) bonds are open with
probability $p$ , while the radial (horizontal) bonds are open with probability $q$. We study the critical curves for these models and establish
 their continuity and strict monotonicity.
\end{abstract}
{\footnotesize Keywords: long range percolation; percolation threshold }

\section{Introduction and Main Result}

Given $k\in\N$, let $\mathbb{S}^k=(\V,\E)$ be the slab of thickness $k$, the graph where the vertex set is
$\V=\Z^2\times\{0,1,...,k\}$ and the set of bonds is $\E=\{\langle x=(x_1,x_2,x_3),y=(y_1,y_2,y_3)\rangle : x,y \in \V, \|x-y\|_1=1\}$, where $\|x-y\|_1=\sum_{i=1}^{3}|x_i-y_i|$
is the usual graph distance in $\Z^3$. The set $\E$ is naturally partitioned in two
disjoint subsets $\E_h$ and $\E_v$. Namely, $\E_h=\{\langle x,y\rangle \in \E:x_3=y_3\}$ and
$\E_v=\{\langle x,y\rangle \in \E:x_1=y_1, x_2=y_2\}$. We say that $e$ is a radial or axial edge acording to $e\in\E_h$ or $\E_v$, respectively.

Given two parameters $p,q\in[0,1]$, we consider a bond anisotropic percolation model on $\mathbb{S}^k$. We associate to each bond $e\in\E$, the state {\em open} or {\em closed} independently, where each bond is open with probability $p$ or $q$, if it belongs to $\E_h$ or $\E_v$, respectively. Thus, this model is described by the probability space
$(\Omega, \mathcal{F},\P_{p,q})$ where $\Omega=\{0,1\}^{\E}$, $\mathcal{F}$ is the
$\sigma$-algebra generated by the cylinder sets in $\Omega$ and $\P_{p,q}=\prod_{e\in\E}\mu(e)$
is the product of Bernoulli measures, where $\mu(e)$ is the Bernoulli measure with parameter $p$
or $q$ acording to $e \in \E_h$ or $e \in \E_v$, respectively. We denote a typical element
of $\Omega$ by $\omega$. When $\omega(e)=1$ we say that $e$ is {\em open}, if $\omega(e)=0$, $e$ is {\em closed}.

Given two vertices $x,y \in \V$ we say that $x$ and $y$ are connected
in the configuration $\omega$ if there exists a finite path of open edges
connecting $x$ to $y$. We will use the short notation $\{x\leftrightarrow y\}$
to denote the set of configurations where $x$ and $y$ are connected.

Given the vertex $x$, the {\em cluster of $x$} in the configuration
$\omega$ is the set $C_x(\omega)=\{y\in \V; x\leftrightarrow y \mbox{ on } \omega\}$. We say that the vertex $x$
percolates when the cardinality of $C_x(\omega)$ is infinite; we will use the following standard notation
$\{x\leftrightarrow\infty\}:=\{\omega\in \Omega;|C_x(\omega)|=\infty\}$, where $|C_x(\omega)|$ is the number
of vertices in $C_x(\omega)$. We define the percolation probability as the function $\theta(p,q):[0,1]^2\mapsto [0,1]$
with $\theta(p,q)=\P_{p,q}(0\leftrightarrow\infty)$. Consider the box $B(n)=[-n,n]^2\times \{0,1,...,k\}$,
denote by $\partial B(n)=\{y=(y_{1},y_{2},y_{3}) \in B(n);\max\{|y_1|,|y_2|\}=n\}$ the boundary of $B(n)$ and $A_n=\{0 \lra \partial B(n)\}$ the event where
$0$ is connected to $\partial B(n)$. Denoting by $\theta_n(p,q)=\P_{p,q}(A_n)$
we observe that $\{0 \lra \infty \}=\cap_{n \geq 1}A_n$ and $\theta_{n}(p,q)\downarrow \theta(p,q)$
as $n$ goes to infinity.

Observe that if $q=1$ the model is equivalent to the bond percolation on $\Z^2$ with parameter $s$ satisfying $1-s=(1-p)^{k+1}$. In this case we denote the horizontal critical value by $p_k$ where $p_k=1-\sqrt[k+1]{\frac{1}{2}}$. Therefore, a simple domination argument shows that $\theta(p,q)=0$ for all $p\leq p_k$ and $q\in[0,1]$.

If $q=0$ we have $k+1$ disjoint copies of $\Z^2$. Then the critical value in this case is $\frac{1}{2}$ and if $p>\frac{1}{2}$ there is percolation, regardless the value of $q$. Using a standard coupling argument we have that $\theta(p,q)$ is non-decreasing function in the parameters $p$ and $q$.
Then we define the function $q^{k}_{c}:[p_k,\frac{1}{2}]\rightarrow [0,1]$ by $$q^{k}_{c}(p)=\sup\{q\in[0,1];\theta(p,q)=0\}.$$

We will show that $q^k_{c}(\frac{1}{2})=0$ and $q^k_{c}(p_k)=1$. As $\theta(p,q)$ is non-decreasing in $p$ then $q^k_{c}(p)$ is non-increasing. In this note we prove that the critical curve $q^k_{c}(p)$ (that divides the regions $\theta(p,q)>0$ and $\theta(p,q)=0$) is continuous and strictly decreasing. Analogously we define the function $p^{k}_{c}:[0,1]\mapsto [p_k,\frac{1}{2}]$ where
\begin{equation}\label{ph}
p^{k}_{c}(q)=\sup\{p\in[0,1];\theta(p,q)=0\}.
\end{equation}

We will omit the index $k$ when it is not necessary and we will write $q_{c}(p)$ and $p_{c}(q)$.

Observe that given some vertex $x\in\V$ the percolation function $\P_{p,q}(x \leftrightarrow\infty)$ is, in general, different from $\theta(p,q)$ but the critical functions $q_{c}(p)$ and $p_{c}(q)$ are the same. Now we can state the main result of this paper:

\begin{teorema}\label{teo1}
The functions $q_{c}:[p_k,\frac{1}{2}]\mapsto [0,1]$  and $p_{c}:[0,1]\mapsto [p_k,\frac{1}{2}]$ are decreasing, continuous and $q_{c}$ is the inverse of $p_{c}$. Moreover, for a compact $[a,b]\subset (p_k,\frac{1}{2})$, there exists positive constants $c(a,b)$ and $C(a,b)$ such that
\begin{equation}\label{lip}
c(a,b)|p'-p|\leq|q_{c}(p')-q_{c}(p)|\leq C(a,b)|p'-p|
\end{equation}
for all $p',p$ $\in [a,b]$.
\end{teorema}

{\bf Remarks:} 1) All the results of this paper can be generalized,
with some minor modifications, for anisotropic percolation in the
whole graph $\Z^3$. 2) Given any $p,p^\prime\in[0,1],\ p\geq
p^\prime$, using the results of Grimmett and Marstrand (see
\cite{GM}, they are still valid for anisotropic percolation), it may be shown that for
all sufficiently large $k$ we have $\theta(p,p^\prime)\geq
\theta(p^\prime,p)$, which says that the slab $\mathbb{S}^k$
percolates better when the greater parameter is on radial bonds. We
expect that this behavior is true for any $k$. Simulations in
\cite{simu} indicates that, in anisotropic $\Z^3$, $q_{c}(p)$ is
convex. If such fact is indeed true for $\mathbb{S}^k$, by Theorem
\ref{teo1}, we have that, if $p>p^\prime$ and $\theta(p,p^\prime)=0$
then $\theta(p^\prime,p)=0$.

In the next section, we state three Lemmas and prove Theorem \ref{teo1}. In Section 3, we will prove the Lemmas stated in Section 2.

\section{Preliminary Lemmas and proof of Theorem 1}

The first lemma proves that for $p<\frac{1}{2}$ and $q$ small enough, we have exponential decay of the radius of the open cluster. As a consequence we have that $q_{c}(p)>0$ for all $p\in [p_k,\frac{1}{2})$.

\begin{lema}\label{dec expo}
Fixed $p<\frac{1}{2}$, there exists $\delta=\delta (p)>0$ with the
following property: For any $0\leq q<\delta$, there exists a
constant $c=c(p,q)>0$ such that
$$\P_{p,q}(|C_{0}|\geq n)\leq e^{-cn},\ \forall n\geq 1.$$ Where
$C_{0}$ is the open cluster of the origin in $S^k$.
\end{lema}

\begin{lema}\label{2}
For all $p\in (p_k,\frac{1}{2}]$, we have that $q_{c}(p)<1$. Moreover, it holds that $\displaystyle\lim_{p\uparrow \frac{1}{2}}q_{c}(p)=$$q_{c}(\frac{1}{2})=0$.
\end{lema}

Combining this result  with Lemma \ref{dec expo}, we have that $0<q_{c}(p)<1$ for all $p\in (p_k,\frac{1}{2})$.

\begin{lema}\label{direcoes}
Given $\delta>0$, there are $\phi=\phi(\delta),\ \psi=\psi(\delta)\in (0,\frac{\pi}{2})$ such that $\forall$ $(p,q)\in [\delta,1-\delta]^2$ the function $\theta(p,q)$ is non-decreasing in the directions of $(\cos\phi,-\sin\phi)$ and $(-\cos\psi,\sin\psi)$.
\end{lema}

 As a consequence of this last lemma, we can see that $\displaystyle\lim_{p\downarrow p_k}q_{c}(p)=q_{c}(p_k)=1$. In fact, suppose that $q_{c}(p_k)<1$ there exists $\bar{q}<1$ such that $\theta(p_k,\bar{q})>0$. Taking $\delta < \min\{\frac{1-\bar{q}}{2},\frac{p_k}{2}\}$, we obtain from Lemma \ref{direcoes} that $\theta(p,q)$ is non-decreasing in the direction $(-\cos\psi,\sin\psi)$. Therefore, it must be a pair $(p,q)\in [\delta,1-\delta]^2$ with $p<p_k$ and $\theta(p,q)>0$, a contradiction, because there is no percolation if $p<p_k$. Then, $q_{c}(p_k)=1$.
 In the same manner we show that $\displaystyle\lim_{p\downarrow p_k}q_{c}(p)=1$.

In resume, we have that $\displaystyle\lim_{p\downarrow p_k}q_{c}(p)=q_{c}(p_k)=1$, $\displaystyle\lim_{p\uparrow \frac{1}{2}}q_{c}(p)=$$q_{c}(\frac{1}{2})=0$ and $0<q_{c}(p)<1,\ \forall  p \in (p_k,\frac{1}{2})$. In analogous way, for the function $p_{c}:[0,1]\mapsto [p_k,1/2]$ defined in Equation \ref{ph}, we have by the facts above that $p_{c}(0)=1/2$, $p_{c}(1)=p_k$ and $p_{c}(q)$ is non-increasing in $q$, then $p_{c}(q)\in [p_k,\frac{1}{2}],\ \forall q \in [0,1]$.

Now, we are able to prove Theorem \ref{teo1}.

\begin{proof}[{\bf Proof of Theorem \ref{teo1}}]

First we will show \eqref{lip}. Fixed $[a,b]\subset (p_k,\frac{1}{2})$ we have $0<q_{c}(b)\leq q_{c}(a)<1$, so take $\delta=\delta(a,b) >0$
such that the square $[2\delta, 1- 2\delta]^2$ contains the points $(a,q_{c}(a))$ and $(b,q_{c}(b))$. As $q_{c}$ is non increasing we have $(p,q_{c}(p)) \in [2\delta,1-2\delta]^2,\ \forall p \in [a,b]$. Let $\phi=\phi(\delta)$ and $\psi=\psi(\delta)$ be given by Lemma \ref{direcoes} and take $\epsilon =\frac{\delta}{2}\min\{(\tan \psi)^{-1},(\tan \phi)^{-1}\}$. Consider then $p<p' \in [a,b]$ with $|p'-p|\leq \epsilon$.

We observe that $\forall$ $0<\eta <\frac{\delta}{2}$

{\it i)} $(p,q_{c}(p)+\eta)$, $(p',q_{c}(p)+\eta-|p'-p|\tan\phi)$
$\in [\delta,1-\delta]^2$ and
\begin{equation}\label{nd}
(p',q_{c}(p)+\eta-|p'-p|\tan\phi)=(p,q_{c}(p)+\eta)+\frac{|p'-p|}{\cos\phi}(\cos\phi,-\sen\phi)
\end{equation}

{\it ii)} $(p,q_{c}(p)-\eta)$, $(p',q_{c}(p)-\eta-|p'-p|\tan\psi)$
$\in [\delta,1-\delta]^2$ and
\begin{equation}\label{ni}
(p',q_{c}(p)-\eta-|p'-p|\tan\psi)=(p,q_{c}(p)-\eta)+\frac{|p'-p|}{\cos\psi}(\cos\psi,-\sen\psi)
\end{equation}

By Lemma \ref{direcoes} and \eqref{nd} we have that
$\theta(p',q_{c}(p)+\eta-|p'-p|\tan\phi)\geq
\theta(p,q_{c}(p)+\eta)>0$, so $q_{c}(p')\leq
q_{c}(p)+\eta-|p'-p|\tan\phi$, $\forall \,0<\eta<\frac{\delta}{2}$,
where we obtain
\begin{equation}
|q_{c}(p')-q_{c}(p)|\geq |p'-p|\tan\phi
\end{equation}

Again, by Lemma \ref{direcoes} (using the consequence that
$\theta(p,q)$ is non-increasing in the direction
$(\cos\psi,-\sen\psi)$) and \eqref{ni}, we have that
$\theta(p',q_{c}(p)-\eta-|p'-p|\tan\psi)\leq
\theta(p,q_{c}(p)-\eta)=0$, so $q_{c}(p')\geq
q_{c}(p)-\eta-|p'-p|\tan\psi,\ \forall 0<\eta<\frac{\delta}{2}$,
where we obtain
\begin{equation}
|q_{c}(p')-q_{c}(p)|\leq |p'-p|\tan\psi
\end{equation}

We proved above the inequalities in \eqref{lip} with the restriction
$|p'-p|\leq\epsilon$. Using that $q_{c}$ is non-increasing we have
that the second inequality in \eqref{lip} is true for all $p',p \in
[a,b]$. So, take $C(a,b)=\tan\psi$. For the first inequality, observe
that if $p'-p>\epsilon$, by taking $p\leq p_1<p_2 \leq p'$ with
$p_2-p_1= \epsilon$ we have
\begin{equation}
|q_{c}(p')-q_{c}(p)|\geq|q_{c}(p_2)-q_{c}(p_1)|\geq
\frac{|p_2-p_1||p'-p|\tan\phi}{|p'-p|}\geq
\frac{\epsilon|p'-p|\tan\phi}{|b-a|}
\end{equation}
so, take $c(a,b)=\frac{\epsilon\tan\phi}{b-a}$.

The argument above shows that $q_{c}$ is Lipschitz-continuous and strictly decreasing on any compact $[a,b] \subset (p_k,\frac{1}{2})$. Combining with $\displaystyle\lim_{p\downarrow p_k}q_{c}(p)=q_{c}(p_k)=1$ and $\displaystyle\lim_{p\uparrow \frac{1}{2}}q_{c}(p)=$$q_{c}(\frac{1}{2})=0$,
we have that $q_{c}$ is strictly decreasing and continuous in the whole interval $[p_k,\frac{1}{2}]$.

Analogously we can prove that $p_{c}(q)$ is strict decreasing and continuous for $q\in [0,1]$.

Now, we will show that the function $q_{c}$ is the inverse function of $p_{c}$. Given any $p\in (p_k,\frac{1}{2})$, as $q_{c}(p)$ is strict decreasing, it holds that $\forall \, \epsilon >0$, $\theta(p-\epsilon,q_{c}(p))=0$ (so $p_{c}(q_{c}(p))\geq p$) and $\theta(p+\epsilon,q_{c}(p))>0$ (so $p_{c}(q_{c}(p))\leq p$). Whence we conclude that $p_{c}(q_{c}(p))=p$. In the same manner we show that $q_{c}(p_{c}(q))=q$. That is, $q_{c}$ is the inverse function of $p_{c}$.

\end{proof}

\section{Proofs of the Lemmas}

\begin{proof} [{\bf  Proof of Lemma \ref{dec expo}}]

We adapt to $\mathbb{S}^k$ the ideas contained in Section 3.5 of \cite{BR}. The key idea for this proof is the Lemma 11 of Section 3.5 in \cite{BR}, which is a consequence of a general result of \cite{LSS} comparing $l$-independent measures with product measures.

Let $G$ be a graph, and let $\tilde{\mathbb{P}}$ be a site percolation measure on $G$, i.e., a probability measure on the set of assignments of states
(open or closed) to the vertices of $G$. The measure $\tilde{\mathbb{P}}$ is $l$-independent if, whenever $U$ and $V$ are sets of vertices of $G$
whose graph distance is at least $l$, the states of the vertices in $U$ are independent of the states of the vertices in $V$. Observe that if $\tilde{\mathbb{P}}$ is $1$-independent, it means that $\tilde{\mathbb{P}}$ is a product measure.

Given an integer $m>0$, denote by $C_m(A)$ the event where some
vertex in $A$ is connected by an open path to a vertex at distance
$m$ from $A$, where the distance is given by the maximum norm. Given
$m>k$ an integer, consider the box $S_m=[0,m-1]^2\times
\{0,1,...,k\}$ in $\mathbb{S}^k$. We will show first that given
$\epsilon>0$ (to be chosen later) and $p<\frac{1}{2}$ there are $m$
and $\delta >0$ such that $\P_{p,q}(C_m(S_m))<\epsilon$, $\forall q
\in [0,\delta)$ .

For all $i\in\{0,1,...,k\}$, we define $P^i$ as the plan
$\Z^2 \times \{i\}$, $S^i_m=[0,m-1]^2\times \{i\} \subset P^i$. Let
$Q(S_m)$ be the event where all vertical bonds of the box
$\tilde{S}_m=[-m,2m-1]^2 \times \{0,...k\}$ are closed.

Observe that $C_m(S_m)\cap Q(S_m)\subset\bigcup_{i=0}^{k}C_m(S^i_m)$. Since
the events $C_m(S^i_m)$ have the same probability, we have
$\P_{p,q}(C_m(S_m)\cap Q(S_m))\leq (k+1) \P_{p,q}(C_m(S^0_m))$.
As $p<\frac{1}{2}$, we use the exponential decay in the
subcritical phase of $\Z^2 \simeq\Z^2\times\{0\}$ (see Theorem 5.4
in \cite{G}). Then, there is some constant $\psi(p)>0$ such that
$\mathbb{P}_{p}(0\leftrightarrow [-n,n]^2)\leq e^{-\psi(p)n}$, where
$\mathbb{P}_{p}$ is the probability measure for ordinary bond
Bernoulli percolation with parameter $p$ on $\Z^2$. Then
\begin{align}
\P_{p,q}(C_m(S^0_m))=\mathbb{P}_{p}(C_m(S^0_m))\leq \displaystyle \sum_{v \in S^0_m}\mathbb{P}_{p}(v \longleftrightarrow \partial B(v,m))\leq \nonumber\\
\leq m^2 e^{-\psi(p)m}\stackrel{m \to \infty}{\longrightarrow} 0
\nonumber
\end{align}
where $B(v,m)$ is the ball of center $v$ and radius $m$ in the
maximum norm in $\Z^2\times \{0\}$. Then, we can take $m$ large
enough such that $\P_{p,q}(C_m(S^0_m))\leq
\frac{\epsilon}{2(k+1)}$. Thereby $\P_{p,q}(C_m(S_m)\cap
Q(S_m))\leq \frac{\epsilon}{2}$.

Observing that $\P_{p,q}(Q(S_m))=(1-q)^N$ where $N=k(3m)^2$ is the
number of vertical bonds in $\tilde{S}_m$, we can choose $\delta>0$
small enough such that $\P_{p,q}(Q(S_m))>1-\frac{\epsilon}{2}$ for
all $p\in[0,1]$ and $q \in [0,\delta)$.

Then fixed $p<\frac{1}{2}$ and $q\in [0,\delta)$, it holds that

$$\P_{p,q}(C_m(S_m))=\P_{p,q}(C_m(S_m)\cap
Q(S_m))+\P_{p,q}(C_m(S_m)\cap Q(S_m)^c)\leq \frac{\epsilon}{2} +
\P_{p,q}(Q(S_m)^c)\leq \epsilon.$$

Now, we will define a site percolation measure $\tilde{\mathbb{P}}$
on $\Z^2$. We declare each vertex $v=(x,y)\in \Z^2$ as open if and
only if the event $C_m(S_{v,m})$ holds for the $m$ by $m$ square
$S_{v,m}=[mx+1,mx+m]\times[my+1,my+m]\times\{0,...,k\}$.

More formally, let $f:\Omega \mapsto \{0,1\}^{\Z^2}$ be the function
defined as $f(\omega)=(f_v(\omega))_{v\in\Z^2}$ where
\begin{equation}
f_v(\omega)=\left\{
\begin{array}
[c]{l}%
1\ ,\mbox{ if}\ \omega\in C_m(S_{v,m})\\
0\ ,\mbox{ if}\ \omega\notin C_m(S_{v,m}).
\end{array}\right.
\end{equation}
The function $f$ and the measure $\P_{p,q}$ induce a probability
measure $\tilde{\mathbb{P}}$ on $\{0,1\}^{\Z^2}$ given by
$\tilde{\mathbb{P}}(A)=\P_{p,q}(f^{-1}(A))$ for any $A\in {\cal
A}$, where ${\cal A}$ is the $\sigma$-algebra generated by the
cylinder sets of $\{0,1\}^{\Z^2}$. This measure $\tilde{\mathbb{P}}$
give us a site percolation model on $\{0,1\}^{\Z^2}$.

Since the event $C_m(S_{v,m})$ depends only on the states of sites
within distance (in the graph distance) $m$ of $S_{v,m}$, then the
measure $\tilde{\mathbb{P}}$ is 5-independent. Furthermore, each
vertex $v \in \Z^2$ is open with probability
$\tilde{\mathbb{P}}(v\mbox{ is open})=\P_{p,q}(C_m(S_{v,m}))\leq
\epsilon$.

From Lemma 11 of Section 3.5 in \cite{BR}, we can take $\epsilon>0$
and $a>0$ (the constants $k$ and $\Delta$ in Lemma 11 are 5 and 4,
respectively) such that if $\tilde{\mathbb{P}}(v\mbox{ is open})\leq
\epsilon$, then $$\tilde{\mathbb{P}}(|\tilde{C}_{v}|\geq n)\leq
e^{-an},\ \forall n\geq 1.$$ Where $\tilde{C}_{v}$ is the open
cluster of the vertex $v$ in this 5-independent model induced on
$\Z^2$.

If $|C_{0}|\geq (k+1)(4m+1)^2$ this implies that every site $u$ of $C_{0}$ is
connected by an open path to some site at distance at least $2m$ from $u$, then $C_m(S_{v,m})$ occurs for every vertex $v\in S_{v,m}\cap C_{0}$, that is, the site $v$ is open in the 5-independent model induced on $\Z^2$, in particular $v\in\tilde{C}_{0}$. Hence, as each $S_{v,m}$ contains $(k+1)m^2$
sites, if $n\geq (k+1)(4m+1)^2$ we have $$\P_{p,q}(|C_{0}|\geq n)\leq
\tilde{\mathbb{P}}\left(|\tilde{C}_{0}|\geq
\frac{n}{(k+1)m^2}\right)\leq e^{-\frac{an}{(k+1)m^2}}.$$

We conclude the proof of Lemma \ref{dec expo} taking $c=\frac{a}{(k+1)m^2}$.

\end{proof}

\begin{proof}[{\bf Proof of Lemma \ref{2}}]

First we will prove that $q_{c}(p)<1$ for $p\in (p_k,1/2]$. Since $q_{c}$ is non-increasing in $p$, it is sufficient to show that $q_{c}(p)<1$ for $p$ close to $p_k$. Given any $\epsilon >0$, we will show that there exists $q<1$ such that que $\theta(p_k+\epsilon,q)>0$, therefore $q_{c}(p_k+\epsilon)\leq q<1$.

Let $u_{1}=(1,0,0),u_{2}=(0,1,0)$ and $u_{3}=(0,0,1)$ be the unitary vectors. Consider the graph $G$ obtained from $\mathbb{S}^k$ replacing each vertical bond $f_{v}=\langle v,v+u_{3}\rangle$, with $v \in \Z^2\times \{0,1,...,k-1\}$, by 4 parallel bonds (denoted by $f_{v}^{r},f_{v}^{l},f_{v}^{t}$ and $f_{v}^{d}$) connecting  the vertices $v$ and $v+u_{3}$ and declaring each of these new bonds open with probability $\widehat{q}$ where $1-q=(1-\widehat{q})^4$.
It means that each vertical bond in $\mathbb{S}^k$ is closed if and only if the respective 4 parallel bonds of $G$ are closed. Observe that $G$ and $\mathbb{S}^k$ have the same vertex set and this replacement does not affect the connective functions involving the vertices of $\mathbb{S}^k$. We have that $\P_{p,q}(0\longleftrightarrow \partial B(n)\mbox{ in }\mathbb{S}^k)=\P_{p,\widehat{q}}(0 \longleftrightarrow \partial B(n)\mbox{ in }G)$ and $\theta^{\mathbb{S}^k}(p,q)=\theta^{G}(p,\widehat{q})$.

We will define a bond percolation process on $\Z^2$ which is stochastically dominated by the bond percolation process on $G$, such that percolation on $\Z^2$ imply percolation on $G$. To simplify the notation we identify $\Z^2$ with $\Z^2 \times \{0\}\subset G$. To each bond $\langle v,v+u_{1}\rangle$ of $\Z^2$ we define the paths $c_{v,u_{1}}^{i}$ on $G$, with $i \in \{0,1,...,k\}$ where
$c_{v,u_{1}}^0$ is the bond $\langle v, v+u_{1}\rangle$ and $c_{v,u_{1}}^i$ for $i=1,\dots,k$ is the path on $G$ that starts at $v$, takes the vertical path using the bonds $f_{v+ju_{3}}^{r}$, with $0\leq j \leq i-1$, until the vertex $v+iu_{3}$, takes the bond $\langle v+iu_{3}, v+iu_{3}+u_{1}\rangle$, and get down vertically until the vertex $v+u_{1}$ using the bonds $f_{v+u_{1}+ju_{3}}^{l}$, with $0\leq j \leq i-1$. Analogously we define the paths $c_{v,u_{2}}^i$ to the bond $\langle v,v+u_{2}\rangle$, using the vertical bonds $f_{v+ju_{3}}^{t}$ and $f_{v+u_{2}+ju_{3}}^{d}$. Declare each bond $e=\langle v, v+u_{1}\rangle$ of $\Z^2$ as open if at least one of the respective paths $c_{v,u_{1}}^i$ is open. Analogously, we do the same for the bond $e=\langle v, v+u_{2}\rangle$ and the paths $c_{v,u_{2}}^i$. Observe that these paths were chosen in such way that we have an independent bond percolation process on $\Z^2$, with parameter $\overline{p}=\overline{p}(p,q)$ defined as:

$$\overline{p}=\overline{p}(p,q)=p\sum_{j=0}^{k}[(1-p)\widehat{q}^2]^j=
p\frac{1-[(1-p)\widehat{q}^2]^{k+1}}{1-(1-p)\widehat{q}^2}$$

reminding that $\widehat{q}=1-(1-q)^\frac{1}{4}$. Taking $p=p_k+\epsilon$ and $q=1$ (so $\widehat{q}=1$) we get $\overline{p}(p_k+\epsilon,1)=1-(\sqrt[k+1]{\frac{1}{2}}-\epsilon)^{k+1}>\frac{1}{2}$. As $\overline{p}(p,q)$ is a continuous function, there is $\delta >0$ such that $\overline{p}(p_k+\epsilon,q)>\frac{1}{2}$, for all $q\in (1-\delta,1]$. Then, we chose $q<1$ such that $\overline{p}(p_k+\epsilon,q)>\frac{1}{2}$. Therefore $\theta^{\mathbb{S}^k}(p,q)=\theta^{G}(p,q)\geq \theta^{\Z^2}(\overline{p}(p_k+\epsilon,q))>0$.

To show that $\displaystyle\lim_{p\uparrow \frac{1}{2}}q_{c}(p)=0$,
we can suppose that $k=1$ since $q^{k}_{c}$ is non-increasing in
$k$. As $\widehat{q}>0$ is equivalent to $q>0$, we have that for
$q>0$ and $p=\frac{1}{2}$,
$\overline{p}=\frac{1}{2}+\frac{\widehat{q}^2}{4}>\frac{1}{2}$.
Fixed $q=\epsilon>0$, as $\overline{p}$ is a continuous function of
$p$, we have that $\overline{p}(p,\epsilon)>\frac{1}{2}$ $\forall \,
p \in (\frac{1}{2}-\frac{\widehat{\epsilon}^2}{4},\frac{1}{2}]$.
Therefore,
$\theta^{\mathbb{S}^k}(p,\epsilon)=\theta^{G}(p,\epsilon)\geq
\theta^{\Z^2}(\overline{p}(p,\epsilon))>0$ and $q_{c}(p)\leq
\epsilon$ $\forall \, p \in
(\frac{1}{2}-\frac{\widehat{\epsilon}^2}{4},\frac{1}{2}]$. We have
then $\displaystyle\lim_{p\uparrow \frac{1}{2}}q_{c}(p)=0$ and
$q_{c}(\frac{1}{2})=0$.
\end{proof}

\begin{proof}[{\bf Proof of Lemma \ref{direcoes}}]

To prove this lemma we enunciated two classical results (without proof) which are adaptations to $\mathbb{S}^k$ of the Russo's formula (see Theorem 2.25 in \cite{G})
and, as consequence, an analogue of Lemma 3.5 of \cite{G}.

Given a configuration $\omega$, we consider the configurations $\omega_e$ and $\omega^{e}$ that coincide with $\omega$ if $f\neq e$,
but $\omega_{e}(e)=0$ and $\omega^{e}(e)=1$ We say that $e$ is pivotal for an increasing event $A$ in the configuration $\omega$ if $I_A(\omega_e)=0$
but $I_A(\omega^e)=1$. We denote by $(e \mbox{ is pivotal for }A)$ the set of such configurations. Thus, as $A_n$ is an increasing event, $e$ is pivotal
for $A_n$ if and only if $A_n$ does not occur when $e$ is closed but $A_n$ does occur when $e$ is open.

\begin{proposicao} \textbf{(Russo's formula)} Consider the function $\theta_{n}(p,q)$, then
\begin{equation}\label{russo}
\begin{array}{l}
\frac{\partial}{\partial p}\theta_{n}(p,q)=\displaystyle\sum_{e \in \mathbb{E}_{h}\cap B(n)}\P_{p,q} (e \mbox{ is pivotal for }A_{n}) \cr \\
\frac{\partial}{\partial q}\theta_{n}(p,q)=\displaystyle\sum_{f \in \mathbb{E}_{v}\cap B(n)}\P_{p,q} (f \mbox{ is pivotal for }A_{n})
\end{array}
\end{equation}
\end{proposicao}

\begin{proposicao}\label{comparacao} There exists a positive integer $N$ and a continuous function $\beta :(0,1)^2\mapsto (0,\infty)$ such that $\forall p,q\in (0,1)$ and $n \geq N$, it holds that
 $$\beta ^{-1}(p,q) \frac{\partial}{\partial p}\theta_{n}(p,q)\geq \frac{\partial}{\partial q}\theta_{n}(p,q) \geq \beta(p,q) \frac{\partial}{\partial p}\theta_{n}(p,q).$$

\end{proposicao}

We show the existence of $\phi$, the proof for $\psi$ is analogue. By the first inequality in Proposition \ref{comparacao}, $\frac{\partial}{\partial p}\theta_{n}(p,q) \geq \beta(p,q) \frac{\partial}{\partial q}\theta_{n}(p,q)$ for all large $n$, where $\beta(p,q)$ is continuous in $(0,1)^2$. Therefore, given $\delta>0$ take $m>0$ such that $\beta(p,q)\geq m$ in $[\delta,1-\delta]^2$, and in this case $\frac{\partial}{\partial p}\theta_{n}(p,q) \geq m \frac{\partial}{\partial q}\theta_{n}(p,q)$. Take then $\phi \in (0,\pi/2)$ such that $\tan\phi=m$. Therefore

\begin{align}\label{dir}
&\nabla \theta_{n}(p,q)\cdot(\cos\phi,-\sin\phi)=\frac{\partial}{\partial p}\theta_{n}(p,q)\cos\phi-
\frac{\partial}{\partial q}\theta_{n}(p,q)\sin\phi=\nonumber\\
&=\cos\phi\left(\frac{\partial}{\partial p}\theta_{n}(p,q)-\frac{\partial}{\partial q}\theta_{n}(p,)\tan\phi\right)\geq \cos\phi\left(m\frac{\partial}{\partial q}\theta_{n}(p,q)-\frac{\partial}{\partial q}\theta_{n}(p,q)\tan\phi\right)=0\nonumber\\
\end{align}
since that $m=\tan\phi$.

Let $(p^\prime,q^\prime)=(p,q)+a(\cos\phi,-\sin\phi)$ such that $(p,q)$ and $(p^\prime,q^\prime)$ are in $[\delta,1-\delta]^2$. Let $\alpha :[0,a]\longrightarrow [\delta,1-\delta]^2$, $\alpha(t)=(p,q)+t(\cos\phi,-\sin\phi)$ the linear path joining $(p,q)$ to $(p^\prime,q^\prime)$. Integrating (\ref{dir}) along the path $\alpha$ we get
$$\theta_{n}(p^\prime,q^\prime)-\theta_{n}(p,q)=\displaystyle\int_{0}^{a}\frac{d}{dt}\theta_{n}(\alpha(t))dt=\displaystyle\int_{0}^{a}\nabla \theta_{n}(p,q)\cdot(\cos\phi,-\sin\phi)dt\geq0.$$
Taking the limit as $n\longrightarrow \infty$ we have $$\theta(p^\prime,q^\prime)=\displaystyle\lim_{n\longrightarrow \infty}\theta_{n}(p^\prime,q^\prime)\geq\lim_{n\longrightarrow \infty}\theta_{n}(p,q)=\theta(p,q)$$

So $\theta(p,q)$ is non-decreasing in the direction $(\cos\phi,-\sin\phi)$.
\end{proof}

{\bf Acknowledgments.} B.N.B. de Lima and R. Sanchis are partially supported by CNPq. R. Sanchis is partially supported by FAPEMIG (Programa Pesquisador Mineiro)

\end{document}